\documentclass[12pt]{article}

\newcommand{\qed}{\hfill$\square$\medskip}

\usepackage{amsmath, epsfig, cite, lineno}
\usepackage{amssymb}
\usepackage{amsfonts}
\usepackage{latexsym}

\newtheorem{thm}{Theorem}[section]

\newtheorem{cor}[thm]{Corollary}

\newtheorem{lem}[thm]{Lemma}

\newtheorem{problem}[thm]{Problem}

\def\pf{\noindent {\it Proof.} }

\numberwithin{equation}{section}

\begin{document}


\begin{center}
{\Large\bf Some congruences related to the $q$-Fermat quotients}
\end{center}

\vskip 2mm \centerline{Victor J. W. Guo}
\begin{center}
{\footnotesize Department of Mathematics, Shanghai Key Laboratory of PMMP,
East China Normal University,\\ 500 Dongchuan Rd., Shanghai 200241,
 People's Republic of China\\
{\tt jwguo@math.ecnu.edu.cn,\quad
http://math.ecnu.edu.cn/\textasciitilde{jwguo} }}
\end{center}

\vskip 0.7cm
{\small \noindent{\bf Abstract.} We give $q$-analogues of the following congruences by Z.-W. Sun:
\begin{align*}
\sum_{k=1}^{p-1}\frac{D_k}{k} \equiv -\frac{2^{p-1}-1}{p} \pmod p,\\
\sum_{k=1}^{p-1}\frac{H_k}{k 2^k}\equiv 0 \pmod{p},\quad p\geqslant 5,
\end{align*}
where $p$ is a prime, $D_n=\sum_{k=0}^{n}{n+k\choose 2k}{2k\choose k}$ are the Delannoy numbers, and
$H_n=\sum_{k=1}^n\frac{1}{k}$ are the harmonic numbers.
We also prove that, for any positive integer $m$ and prime $p>m+1$,
\begin{align*}
\sum_{1\leqslant k_1\leqslant \cdots \leqslant k_m\leqslant p-1}\frac{1}{k_1\cdots k_m 2^{k_m}}
\equiv\frac{1}{2}\sum_{k=1}^{p-1}\frac{(-1)^{k-1}}{k^m} \pmod p,
\end{align*}
which is a multiple generalization of Kohnen's congruence. Furthermore, a $q$-analogue of this congruence is
established.                                   }


\vskip 5mm
\noindent{\it MR Subject Classifications}: 11B65; 33D15; 05A10; 05A30

\section{Introduction}
Fermat's Little Theorem states that if $p$ is a prime, then for any integer $a$ not divisible by $p$,
the number $a^{p-1}-1$ is a multiple of $p$.
Numbers of the form $(a^{p-1}-1)/p$ are called Fermat quotients of $p$ to base $a$.
There are several different congruences for the Fermat quotients $(2^{p-1}-1)/p$ in the literature.

Let $p\geqslant 3$ be a prime. Since $\frac{1}{p}{p\choose k}\equiv \frac{(-1)^{k-1}}{k} \pmod{p}$
for $1\leqslant k\leqslant p-1$, the following result is well-known:
\begin{align}\label{eq:known}
\frac{2^{p-1}-1}{p}\equiv \frac{1}{2}\sum_{k=1}^{p-1}\frac{(-1)^{k-1}}{k}
\equiv \sum_{k=1}^{(p-1)/2}\frac{(-1)^{k-1}}{k} \pmod{p}.
\end{align}
 A classical Glaisher's congruence (see \cite{Glaisher,Granville}) for Fermat quotients is
\begin{align}\label{eq:Glaisher}
\sum_{k=1}^{p-1}\frac{2^{k-1}}{k}
\equiv -\frac{2^{p-1}-1}{p}\pmod p.
\end{align}
Kohnen \cite{Kohnen} established the following congruence
\begin{align}\label{eq:Kohnen}
\sum_{k=1}^{p-1}\frac{1}{k 2^k}
\equiv \sum_{k=1}^{(p-1)/2}\frac{(-1)^{k-1}}{k}\pmod p.
\end{align}
Z.-W. Sun \cite{Sun} prove that
\begin{align}\label{eq:Sun}
\sum_{k=1}^{p-1}\frac{D_k}{k} \equiv -\frac{2^{p-1}-1}{p} \pmod p,
\end{align}
where $D_n$ are the (central) Delannoy numbers defined by
$$
D_n=\sum_{k=0}^{n}{n+k\choose 2k}{2k\choose k}.
$$

Pan \cite{Pan} gave $q$-analogues of \eqref{eq:known} and \eqref{eq:Glaisher} as follows:
\begin{align}
\sum_{k=1}^{p-1}\frac{(-1)^k}{[k]}
&\equiv -\frac{2(-q;q)_{p-1}-2}{[p]}-\frac{(p-1)(1-q)}{2}\pmod{[p]}, \label{eq:cong-2}  \\
\sum_{k=1}^{p-1} \frac{(-q;q)_{k}q^k}{2[k]}
&\equiv -\frac{(-q;q)_{p-1}-1}{[p]}-\frac{(p-1)(1-q)}{2} \pmod{[p]}, \label{eq:pan-1}
\end{align}
where $(a;q)_n=(1-a)(1-aq)\cdots(1-aq^{n-1})$ and $[n]=(1-q^n)/(1-q)$. The polynomials
$$
\frac{(-q;q)_{p-1}-1}{[p]}
$$
are called the $q$-Fermat quotients of $p$ to base $2$.
Tauraso \cite{Tauraso}
obtained the following $q$-analogues of \eqref{eq:Glaisher} and \eqref{eq:Kohnen}:
\begin{align*}
\sum_{k=1}^{p-1} \frac{(-q;q)_{k-1}q^{-{k\choose 2}}}{[k]}
&\equiv -\frac{(-q;q)_{p-1}-1}{[p]} \pmod{[p]},    \\
\sum_{k=1}^{p-1}\frac{q^k}{[k](-q;q)_k}
&\equiv \frac{(-q;q)_{p-1}-1}{[p]} \pmod{[p]}.  \notag
\end{align*}

In this paper, we give a $q$-analogue of \eqref{eq:Sun}
and new $q$-analogues of \eqref{eq:Glaisher} and \eqref{eq:Kohnen}.
\begin{thm}\label{thm:main1}
For any prime $p\geqslant 3$, there holds
\begin{align}
\sum_{k=1}^{p-1} \frac{D_k(q)}{[k]}
\equiv -\frac{(-q;q)_{p-1}-1}{[p]}+\frac{(p-1)(1-q)}{4} \pmod{[p]}, \label{eq:Delannoy}
\end{align}
Here the $q$-Delannoy numbers $D_n(q)$ are defined  by
$$
D_n(q)=\sum_{k=0}^n\frac{1+q^k}{2}{n+k\brack 2k}{2k\brack k}q^{{k\choose 2}-2nk},
$$
where
$$
{n\brack k}=\frac{(q;q)_n}{(q;q)_k (q;q)_{n-k}}
$$
stands for the $q$-binomial coefficient.
\end{thm}

\begin{thm}\label{thm:main2}
For any prime $p\geqslant 3$, there hold
\begin{align}
\sum_{k=1}^{p-1}\frac{(-q;q)_{k-1} q^k}{[k]}
&\equiv -\frac{(-q;q)_{p-1}-1}{[p]}-\frac{(p-1)(1-q)}{2}
\pmod{[p]},  \label{eq:first-p} \\
\sum_{k=1}^{p-1}\frac{q^{k+1\choose 2}}{[k](-q;q)_k}
&\equiv \frac{(-q;q)_{p-1}-1}{[p]}+\frac{(p-1)(1-q)}{2} \pmod{[p]}. \label{eq:cong-K}
\end{align}
\end{thm}

We shall also give a multiple generalization of \eqref{eq:Kohnen} as follows.
\begin{thm}\label{thm:main3}
For any positive integer $m$ and prime $p> m+1$, there holds
\begin{align}
\sum_{1\leqslant k_1\leqslant \cdots \leqslant k_m\leqslant p-1}\frac{1}{k_1\cdots k_m 2^{k_m}}
\equiv\frac{1}{2}\sum_{k=1}^{p-1}\frac{(-1)^{k-1}}{k^m} \pmod p.  \label{eq:multi-cong}
\end{align}
In particular, if $m$ is even, then
\begin{align}
\sum_{1\leqslant k_1\leqslant \cdots \leqslant k_m\leqslant p-1}\frac{1}{k_1\cdots k_m 2^{k_m}}
\equiv 0\pmod p. \label{eq:multi-even}
\end{align}
\end{thm}
Note that, when $m=2$, the congruence \eqref{eq:multi-even} can be written as
\begin{align}
\sum_{k=1}^{p-1}\frac{H_k}{k 2^k}\equiv 0 \pmod{p},\quad p\geqslant 5, \label{eq:sun-harmonic}
\end{align}
where $H_n=1+\frac{1}{2}+\cdots+\frac{1}{n}$ are the harmonic numbers. The congruence \eqref{eq:sun-harmonic}
was first proved by Z.-W. Sun \cite{Sun} and generalized to the modulus $p^2$ case by Sun and Zhao \cite{SZ}.
Some other generalizations and refinements of \eqref{eq:sun-harmonic} can be found in \cite{Sun12-2,Me0}.

Let
$$H_n(q)=\sum_{k=1}^n\frac{1}{[k]}$$
be the $q$-harmonic numbers.
Our last theorem is the following neat $q$-analogue of \eqref{eq:sun-harmonic}.
\begin{thm}\label{thm:main4}
For any prime $p\geqslant 5$, there holds
\begin{align*}
\sum_{k=1}^{p-1}\frac{H_k(q)q^{k+1\choose 2}}{[k](-q;q)_k }
\equiv  \frac{(p^2-1)(1-q)^2}{24} \pmod{[p]}.
\end{align*}
\end{thm}

The paper is organized as follows. In the next section, we give a proof of Theorem \ref{thm:main1}
by using some $q$-series identities and known $q$-congruences.
In Section 3, we give proofs of Theorems \ref{thm:main2}--\ref{thm:main4}
by first establishing a multiple series generalization of Kohnen's identity \cite{Kohnen}:
\begin{align}
\sum_{k=1}^{n}\frac{1}{k}(1-x)^k
=\sum_{k=1}^n\frac{(-1)^k}{k}{n\choose k}(x^k-1).  \label{eq:Kbino}
\end{align}
In fact, a $q$-analogue of \eqref{eq:multi-cong} will be proved.
Some consequences and remarks will be mentioned in the last section.

\section{Proof of Theorem \ref{thm:main1}}
Applying the Lagrange interpolation formula for $x^r$ at the values $q^{-k}$ ($0\leqslant k\leqslant n$)
of $x$, we have the following result (see \cite[Theorem 1.1]{GZpq} for a generalization),
which will be used in the proof of Theorem \ref{thm:main1}.
\begin{lem}\label{lem:lag}
For $n\geqslant 1$ and $0\leqslant r\leqslant n$, there holds
\begin{align}
\sum_{k=0}^n (-1)^k {n\brack k} \frac{q^{{k+1\choose 2}-rk}}{1-xq^k}
=\frac{(q;q)_n}{(x;q)_{n+1}}x^r. \label{eq:lag}
\end{align}
\end{lem}

\medskip
\noindent{\it Proof of Theorem {\rm\ref{thm:main1}}.}
By the $q$-Lucas theorem (see \cite{De, Olive, GZ}), or by the factorization
of $q$-binomial coefficients into cyclotomic polynomials (see \cite{CH,KW}),
for any prime $p\geqslant 3$ and $(p-1)/2<k<p$, there holds
$$
{2k\brack k}\equiv 0 \pmod{[p]}.
$$
Hence, by the $q$-Chu-Vandermonde identity (see \cite[(3.3.10)]{Andrews}), we have
\begin{align*}
\sum_{m=1}^{p-1} \frac{D_m(q)-1}{[m]}
&=\sum_{m=1}^{p-1}\sum_{k=1}^{m}\frac{1+q^k}{2[m]}{2k\brack k}{m+k\brack 2k}q^{{k\choose 2}-2mk} \\
&=\sum_{k=1}^{p-1}\frac{1+q^k}{2}{2k\brack k}q^{k\choose 2}\sum_{m=k}^{p-1}
\sum_{j=k}^{2k}\frac{q^{-j(2k-j+m)}}{[m]}{m\brack j}{k\brack 2k-j} \\
&\equiv\sum_{k=1}^{(p-1)/2}\frac{1+q^k}{2}{2k\brack k}q^{k\choose 2}
\sum_{m=k}^{p-1}\sum_{i=k}^{2k}\frac{q^{-j(2k-j+m)}}{[m]}{m\brack j}{k\brack 2k-j} \pmod{[p]}.
\end{align*}
Note that
\begin{align}
\frac{1}{[m]}{m\brack j}&=\frac{1}{[j]}{m-1\brack j-1},  \notag\\
\sum_{m=j}^{p-1}{m-1\brack j-1}q^{-mj} &={p-1\brack j}q^{-(p-1)j}, \notag\\
{p-1\brack k}=\prod_{j=1}^{k}\frac{1-q^{p-j}}{1-q^j}
&\equiv \prod_{j=1}^{k}\frac{1-q^{-j}}{1-q^j}=(-1)^k q^{-{k+1\choose 2}} \pmod{[p]}. \label{eq:qbino-cong}
\end{align}
For $1\leqslant k\leqslant (p-1)/2$, we have
\begin{align*}
&\hskip -3mm\sum_{m=k}^{p-1}\sum_{j=k}^{2k}\frac{q^{-j(2k-j+m)}}{[m]}{m\brack j}{k\brack 2k-j} \\
&=\sum_{j=k}^{2k}\frac{q^{-j(2k-j)}}{[j]}{k\brack j-k}
\sum_{m=j}^{p-1}{m-1\brack j-1}q^{-mj} \\
&\equiv \sum_{j=k}^{2k}\frac{q^{-j(2k-j+p-1)}}{[j]}{k\brack j-k}{p-1\brack j}\\
&\equiv
\sum_{i=k}^{2k}(-1)^j\frac{q^{{j+1\choose 2}-2jk}}{[j]}{k\brack j-k} \pmod{[p]}.
\end{align*}
By Lemma \ref{lem:lag}, we have
\begin{align*}
\sum_{j=k}^{2k}(-1)^j\frac{q^{{j+1\choose 2}-2jk}}{[j]}{k\brack j-k}
&=(1-q)q^{-k(3k-1)/2}\sum_{j=k}^{2k}(-1)^j\frac{q^{{j-k+1\choose 2}-(j-k)k}}{1-q^{(j-k)+k}}{k\brack j-k} \\
&=\frac{(-1)^k(1-q)(q;q)_k q^{-k(3k-1)/2+k^2}}{(q^k;q)_{k+1}} \\
&=\frac{(-1)^k q^{-{k\choose 2}}}{[k]{2k\brack k}}.
\end{align*}
It follows that
\begin{align*}
\sum_{m=1}^{p-1} \frac{D_m(q)-1}{[m]}
&\equiv \sum_{k=1}^{(p-1)/2}\frac{(-1)^k(1+q^k)}{2[k]}  \\
&\equiv \frac{1}{2}\sum_{k=1}^{(p-1)/2}\frac{(-1)^k}{[k]}
+\frac{1}{2}\sum_{k=1}^{(p-1)/2}\frac{(-1)^{p-k}}{[p-k]} \\
&= \frac{1}{2}\sum_{k=1}^{p}\frac{(-1)^k}{[k]}\pmod{[p]}.
\end{align*}
The proof then follows from \eqref{eq:cong-2} and the following congruence due to Andrews \cite{Andrews99}:
\begin{align}
H_{p-1}(q)
\equiv \frac{(p-1)(1-q)}{2} \pmod{[p]}. \label{eq:Andrews}
\end{align}
\qed

\noindent{\it Remark.} If we define the $q$-Delannoy numbers by
$$
\overline{D}_n(q)=\sum_{k=0}^n {n+k\brack 2k}{2k\brack k}q^{{k\choose 2}-2nk},
$$
then we have the following congruence:
\begin{align}
\sum_{m=1}^{p-1} \frac{\overline{D}_m(q)-1}{[m]}
\equiv \sum_{k=1}^{(p-1)/2}\frac{(-1)^k}{[k]} \pmod{[p]}.  \label{eq:rmk}
\end{align}
However, it is difficult to determine the right-hand side of \eqref{eq:rmk} modulo $[p]$.
This is why we need to replace $\overline{D}_n(q)$ by $D_n(q)$ in Theorem \ref{thm:main1}.

\section{Proofs of Theorems \ref{thm:main2}--\ref{thm:main4}}
Dilcher \cite{Dilcher} established the following identity:
\begin{align}
\sum_{1\leqslant k_1\leqslant\cdots\leqslant k_m\leqslant n} \frac{q^{k_1+\cdots+k_m}}{(1-q^{k_1})\cdots(1-q^{k_m})}
=\sum_{k=1}^{n}(-1)^{k-1}{n\brack k}\frac{q^{{k\choose 2}+km}}{(1-q^k)^m}, \label{eq:dilch}
\end{align}
which is a multiple series generalization of  Van Hamme's identity \cite{Hamme}:
\begin{align}
\sum_{k=1}^n \frac{q^k}{1-q^k}
=\sum_{k=1}^{n}(-1)^{k-1}{n\brack k}_q\frac{q^{k+1\choose 2}}{1-q^k}, \label{eq:hamme}
\end{align}
Further generalizations of Dilcher's identity \eqref{eq:dilch} have been obtained by Fu and Lascoux
\cite{FL,FL2}, Zeng \cite{Zeng}, Ismail and Stanton \cite{IS}, Guo and Zhang \cite{GZh}, Gu and Prodinger \cite{GP},
and Guo and Zeng \cite{GZ}.

In what follows we give a new generalization of Dilcher's identity \eqref{eq:dilch}
that also include Kohnen's identity \eqref{eq:Kbino} as a special case.
\begin{thm}For $m,n\geqslant 1$, there holds
\begin{align}
\sum_{1\leqslant k_1\leqslant\cdots\leqslant k_m\leqslant n}
\frac{(x;q)_{k_1}q^{k_1+\cdots+k_m}}{(1-q^{k_1})\cdots(1-q^{k_m})}
=\sum_{k=1}^{n}(-1)^k {n\brack k} \frac{q^{{k\choose 2}+km}}{(1-q^k)^m}  (x^k-1). \label{eq:xDilcher}
\end{align}
\end{thm}
\pf For $1\leqslant r\leqslant n$, the coefficient of $x^r$ in the left-hand side of \eqref{eq:xDilcher}
is given by
\begin{align}
(-1)^r q^{r\choose 2}\sum_{k_m=r}^n \sum_{k_{m-1}=r}^{k_m}\cdots \sum_{k_2=r}^{k_3}\sum_{k_1=r}^{k_2}
{k_1\brack r}\frac{q^{k_1+\cdots+k_m}}{(1-q^{k_1})\cdots(1-q^{k_m})}. \label{eq:multi}
\end{align}
It is easy to see that
\begin{align}
\sum_{k_1=r}^{k_2}{k_1\brack r}\frac{q^{k_1}}{1-q^{k_1}}
=\frac{1}{1-q^r}\sum_{k_1=r}^{k_2}{k_1-1\brack r-1}q^{k_1}
=\frac{1}{1-q^r}{k_2\brack r}q^{r}.  \label{eq:k1tok2}
\end{align}
By repeatedly using the summation formula \eqref{eq:k1tok2}, one sees that \eqref{eq:multi}
is equal to
\begin{align}
(-1)^r  \frac{q^{{r\choose 2}+mr}}{(1-q^r)^m} {n\brack r}. \label{eq:multi-sum}
\end{align}
That is to say, the coefficients of $x^r$ in both sides of \eqref{eq:xDilcher}
are equal for $1\leqslant r\leqslant n$.
Also \eqref{eq:xDilcher} is true for $x=1$. Therefore it must be true for any $x$.
\qed

\noindent{\it Remark.} An equivalent form of the fact that \eqref{eq:multi} equals \eqref{eq:multi-sum}
has been given by Fu and Lascoux \cite[Lemma 2.1]{FL2}. The proof given here is
more straightforward.
\medskip

When $m=1$, we obtain the following result, which is crucial in the proof of Theorem \ref{thm:main2}.
\begin{cor}For $n\geqslant 1$, there holds
\begin{align}
\sum_{k=1}^{n}\frac{(x;q)_k q^k}{1-q^k} =\sum_{k=1}^{n}\frac{(-1)^k
q^{k+1\choose 2}}{1-q^k}  {n\brack k}(x^k-1). \label{eq:qKohnen}
\end{align}
\end{cor}

\noindent{\it Proof of Theorem {\rm\ref{thm:main2}}.}
In \eqref{eq:qKohnen} we set $n=p-1$ and $x=-1$ and multiply both sides by $1-q$.
By \eqref{eq:qbino-cong} and \eqref{eq:Andrews},
the equation \eqref{eq:qKohnen} simplifies to
\begin{align}
\sum_{k=1}^{p-1}\frac{(-1;q)_k q^k}{[k]}
\equiv\sum_{k=1}^{p-1}\frac{(-1)^k }{[k]}-\frac{(p-1)(1-q)}{2}
\pmod{[p]}, \label{eq:qK-p-1}
\end{align}
which is the desired congruence \eqref{eq:first-p}.
Replacing $k$ by $p-k$ on the left-hand side of \eqref{eq:qK-p-1}, we have
\begin{align*}
\sum_{k=1}^{p-1}\frac{(-1;q)_k q^k}{[k]}
&\equiv -2\sum_{k=1}^{p-1}\frac{(-q;q)_{p-k-1}}{[k]}  \\
&= -2\sum_{k=1}^{p-1}\frac{(-q;q)_{p-1}}{[k](-q^{p-k};q)_k} \\
&\equiv -2\sum_{k=1}^{p-1}\frac{q^{k+1\choose 2}}{[k](-q;q)_k}  \pmod{[p]},
\end{align*}
By \eqref{eq:cong-2}, the right-hand side of \eqref{eq:qK-p-1} is congruent to
$$
-\frac{(-q;q)_{p-1}-1}{[p]}-\frac{(p-1)(1-q)}{2} \pmod{[p]}.
$$
This proves \eqref{eq:cong-K}.  \qed

Noticing that
\begin{align*}
\sum_{k=1}^{p}\frac{(-1)^k}{[k]}
&=\sum_{k=1}^{(p-1)/2}\frac{(-1)^k}{[k]}
+\sum_{k=1}^{(p-1)/2}\frac{(-1)^{p-k}}{[p-k]}
\equiv \sum_{k=1}^{(p-1)/2}\frac{(-1)^k(1+q^k)}{[k]}
\pmod{[p]},
\end{align*}
we can also rewrite \eqref{eq:cong-K} as
$$
\sum_{k=1}^{p-1}\frac{q^{k+1\choose 2}}{[k](-q;q)_k}
\equiv\sum_{k=1}^{(p-1)/2}\frac{(-1)^{k-1}(1+q^k)}{2[k]}+\frac{(p-1)(1-q)}{4} \pmod{[p]}.
$$

Similarly, if we set $n=p-1$ and $x=-1$, multiply both sides by $(1-q)^m$ in \eqref{eq:xDilcher},
replace $k_i$ by $p-k_i$ for $1\leqslant i\leqslant m$, and finally reverse the order of
$k_1,\ldots,k_m$, then we obtain the following result.

\begin{thm}
For any positive integer $m$ and prime $p\geqslant 3$, there holds
\begin{align}
\sum_{1\leqslant k_1\leqslant \cdots\leqslant k_m\leqslant p-1}
\frac{q^{{k_m+1\choose 2}}}{[k_1]\cdots [k_m](-q;q)_{k_m}}
\equiv (-1)^{m} \sum_{k=1}^{p-1} \frac{q^{(m-1)k}}{2[k]^m}  ((-1)^k-1)   \pmod{[p]}.   \label{eq:xxxx}
\end{align}
\end{thm}
\noindent{\it Proof of Theorem {\rm\ref{thm:main3}}}.
Letting $q=1$ in \eqref{eq:xxxx} and
using the classical congruence
$$
\sum_{k=1}^{p-1}\frac{1}{k^m}\equiv 0 \pmod{p}, \quad\text{for $p>m+1$},
$$
we complete the proof of \eqref{eq:multi-cong}.

For $m$ even, replacing $k$ by $p-k$, one sees that
\begin{align*}
\sum_{k=1}^{p-1}\frac{(-1)^{k}}{k^m}\equiv \sum_{k=1}^{p-1}\frac{(-1)^{p-k}}{(p-k)^m}
\equiv -\sum_{k=1}^{p-1}\frac{(-1)^{k}}{k^m} \equiv 0 \pmod{p}.
\end{align*}
This proves \eqref{eq:multi-even}.
\qed

\noindent{\it Proof of Theorem {\rm\ref{thm:main4}}}.
When $m=2$, the congruence \eqref{eq:xxxx} can be written as
\begin{align}
\sum_{k=1}^{p-1}\frac{H_k(q)q^{k+1\choose 2}}{[k](-q;q)_k }
\equiv  \sum_{k=1}^{p-1} \frac{(-1)^{k} q^{k}}{2[k]^2}
-\sum_{k=1}^{p-1} \frac{q^{k}}{2[k]^2}\pmod{[p]},\quad\text{for } p\geqslant 5. \label{eq:sun-harmonic-q}
\end{align}
Note that
\begin{align}
\sum_{k=1}^{p-1}\frac{(-1)^{k}q^k}{[k]^2}\equiv \sum_{k=1}^{p-1}\frac{(-1)^{p-k}q^{p-k}}{[p-k]^2}
\equiv -\sum_{k=1}^{p-1}\frac{(-1)^{k}q^k}{[k]^2} \equiv 0 \pmod{[p]}, \label{eq:q-zero}
\end{align}
and Shi and Pan \cite[(4)]{SP} proved that
\begin{align}
\sum_{k=1}^{p-1}\frac{q^k}{[k]^2}\equiv -\frac{(p^2-1)(1-q)^2}{12} \pmod{[p]}. \label{eq:Shi-Pan}
\end{align}
The proof then follows form combing \eqref{eq:sun-harmonic-q}--\eqref{eq:Shi-Pan}.
\qed

\section{Some consequences and remarks}

\begin{cor}For any prime $p\geqslant 3$, there holds
\begin{align}
\sum_{k=1}^{p-1} \frac{kq^k}{1+q^k}\equiv \frac{p(p-1)(1-q)}{2}+p\sum_{k=1}^{p-1} \frac{(-q;q)_{k-1}q^k}{1-q^k}
\pmod{[p]}. \label{eq:-thm2}
\end{align}
\end{cor}
\pf Multiplying both sides of \eqref{eq:first-p} by $1-q^p$, we have
\begin{align}
(-q;q)_{p-1}-1\equiv
-(1-q^p)\left(\frac{(p-1)(1-q)}{2}+\sum_{k=1}^{p-1} \frac{(-q;q)_{k-1}q^k}{[k]}\right)
\pmod{[p]^2}. \label{eq:second-p}
\end{align}
Differentiating both sides of \eqref{eq:second-p} with respect to $q$, we obtain
\begin{align}
(-q;q)_{p-1}\sum_{k=1}^{p-1} \frac{kq^{k-1}}{1+q^k} \equiv
pq^{p-1}\left(\frac{(p-1)(1-q)}{2}+\sum_{k=1}^{p-1} \frac{(-q;q)_{k-1} q^k}{[k]}\right). \label{eq:third-p}
\end{align}
Since $(-q;q)_{p-1}\equiv q^{p}\equiv 1\pmod{[p]}$, one sees that \eqref{eq:third-p} is
equivalent to \eqref{eq:-thm2}.
\qed

Combining \eqref{eq:first-p} and \eqref{eq:-thm2}, we obtain the following result.
\begin{cor}For any prime $p\geqslant 3$, there holds
\begin{align*}
\frac{(-q;q)_{p-1}-1}{1-q^p}\equiv
-\frac{1}{p}\sum_{k=1}^{p-1} \frac{kq^k}{1+q^k} \pmod{[p]}.
\end{align*}
\end{cor}

Letting $q\to 1$, $n=p-1$, and $x\in\mathbb{Z}$ in \eqref{eq:xDilcher}, we get
\begin{cor}
For any integer $x$, positive integer $m$, and prime $p\geqslant 3$, there holds
\begin{align}
\sum_{1\leqslant k_1\leqslant\cdots\leqslant k_m\leqslant p-1}
\frac{(1-x)^{k_1}}{k_1\cdots k_m}
\equiv \sum_{k=1}^{p-1} \frac{(x^k-1)}{k ^m} \pmod{p}. \label{eq:xxyy}
\end{align}
\end{cor}
Letting $x=-1$ or $x=2$ in \eqref{eq:xxyy}, we have
\begin{align}
\sum_{1\leqslant k_1\leqslant\cdots\leqslant k_m\leqslant p-1}
\frac{2^{k_1}}{k_1\cdots k_m}
&\equiv \sum_{k=1}^{p-1} \frac{((-1)^k-1)}{k ^m} \pmod{p}, \label{eq:cor1} \\
\sum_{1\leqslant k_1\leqslant\cdots\leqslant k_m\leqslant p-1}
\frac{(-1)^{k_1}}{k_1\cdots k_m}
&\equiv \sum_{k=1}^{p-1} \frac{(2^k-1)}{k ^m} \pmod{p}  \label{eq:cor2}
\end{align}
It follows from \eqref{eq:cor1} and \eqref{eq:cor2} that
\begin{align}
\sum_{1\leqslant k_1\leqslant\cdots\leqslant k_m\leqslant p-1}
\frac{2^{k_1}-(-1)^{k_1}}{k_1\cdots k_m}
\equiv \sum_{k=1}^{p-1} \frac{((-1)^k-2^k)}{k ^m} \pmod{p}. \label{eq:cor3}
\end{align}
The $m=3$ case of \eqref{eq:cor3} has already appeared in \cite{SZ}.

Z.-W. Sun \cite{Sun95} proved that
\begin{align}\label{eq:Sun95}
\sum_{k=1}^{(p-1)/2}\frac{1}{k 2^k}
\equiv \sum_{k=1}^{\lfloor 3p/4\rfloor}\frac{(-1)^{k-1}}{k}\pmod p,
\end{align}
where $\lfloor x\rfloor$ denotes the greatest integer not exceeding $x$. On the other hand, Z.-H. Sun \cite[Theorem 4.1(iii)]{SunZH} gave a generalization of Kohnen's congruence \eqref{eq:Kohnen}
as follows:
\begin{align}\label{eq:SunZH}
\sum_{k=1}^{p-1}\frac{1}{k 2^k}
\equiv \frac{2^{p-1}-1}{p}-\frac{(2^{p-1}-1)^2}{2p}\pmod{p^2}.
\end{align}
Recently, an elementary proof of \eqref{eq:SunZH} has been given by Me\v{s}trovi\'{c} \cite{Me}.

We end the paper with the following problem.
\begin{problem}
Are there $q$-analogues of the congruences \eqref{eq:Sun95} and \eqref{eq:SunZH}?
\end{problem}

\vskip 3mm
\noindent{\bf Acknowledgment.} This work was partially
supported by the Fundamental Research Funds for the Central Universities and
the National Natural Science Foundation of China (grant 11371144).

\renewcommand{\baselinestretch}{1}

\end{document}